\documentclass{commat}

\newcommand{\R}{{\mathbb{R}}}
\newcommand{\Z}{{\mathbb{Z}}}

\title{%
An existence result for $p$-Laplace equation with gradient nonlinearity in $\R^N$
}

\author{%
Shilpa Gupta and Gaurav Dwivedi
}

\affiliation{
    \address{%
   Department of Mathematics. Birla Institute of Technology and Science Pilani. Pilani Campus, Vidya Vihar. Pilani, Jhunjhunu. Rajasthan, India - 333031.
        }
    \email{%
    p20180442 and gaurav.dwivedi @pilani.bits-pilani.ac.in
    }
    }

\abstract{
We prove the existence of a weak solution to the problem 
\begin{equation*}
\begin{split}
-\Delta_{p}u+V(x)|u|^{p-2}u & =f(u,|\nabla u|^{p-2}\nabla u),  \ \ \ \\
u(x) & >0\ \   \forall x\in\R^{N},
\end{split}
\end{equation*}
where $\Delta_{p}u=\hbox{div}(|\nabla u|^{p-2}\nabla u)$ is the $p$-Laplace operator, $1<p<N$ and the nonlinearity  $f:\R\times\R^{N}\rightarrow\R$ is  continuous and it depends on gradient of the solution. We use an iterative technique based on the Mountain pass theorem to prove our existence~result.
    }

    

\VOLUME{30}

\NUMBER{1}

\firstpage{149}

\DOI{https://doi.org/10.46298/cm.9316}

\begin{paper}

\section{Introduction}
\setcounter{equation}{0}
In this article, we  prove the existence of a weak solution to the problem:
\begin{gather} \label{1.1}
-\Delta_{p}u+V(x)|u|^{p-2}u = f(u,|\nabla u|^{p-2}\nabla u), \notag \\
u(x) > 0 \quad \forall x\in\R^{N},
\end{gather}
where $1<p<N$ and the non-linearity $f:\R\times\R^{N}\rightarrow\R$ is a continuous function. 

 The  Problem (\ref{1.1}) is non-variational in nature,  as the  nonlinearity  $f$  depends on gradient of the solution. Such type of problems have been studied widely in literature through non-variational techniques, such as method of sub-solution and super-solution \cite{Faria}, \cite{serva}, \cite{xavier}, degree theory \cite{ruiz}, \cite{wang} etc. In 2004,  Figueiredo $et \ al.$ \cite{figu}  used an iterative technique based on Mountain Pass Theorem to establish the existence of a positive and a negative solution to the problem:
\begin{equation*}
\begin{split}
-\Delta u &=f(x,u,\nabla u)\,\,\,\,\text{ in }\Omega,   \\
u&=0\ \,\,\,\,\, \text{ on }\partial \Omega,
\end{split}
\end{equation*}
where $\Omega\subseteq \R^n$ is a  smooth and bounded domain. Motivated by the techniques used by Figueiredo et. al \cite{figu}, several authors established existence results for second order elliptic equations  with gradient nonlinearities, see for instance \cite{chung}, \cite{gfigu}, \cite{girardi}, \cite{liu2017}, \cite{pimenta}, \cite{serva2} and references therein. 

This work is motivated the existence results of G.M. Figueiredo \cite{gfigu}, where the author obtained existence of a positive solution to \eqref{1.1} with $V(x)=1.$ Recently, an existence result for \eqref{1.1} in case of $p=N$ is discussed by Chen et al. \cite{chen}. For some existence results for the problems of the type \eqref{1.1} with potential $V(x)$ and without gradient dependence, we refer to \cite{alves}, \cite{embed}, \cite{evaraldo} and references therein.  

The plan of this article is as follows: In section 2, we state our hypotheses and main result. Section 3 deals with the proof of our main result, i.e., Theorem \ref{t1}.

\section{Hypotheses and Main Result}
In this section, we state hypotheses on the nonlinearity $f$ and the potential $V.$
We assume the following conditions on the nonlinearity  $f$:
\begin{itemize}
\item[$ (f_{1}) $] $f(t,|\xi|^{p-2}\xi)=0$ for all $t<0, \ \xi\in\R^{N}.$
\item[$ (f_{2}) $]  $\displaystyle\lim_{|t|\rightarrow 0}\dfrac{|f(t,|\xi|^{p-2}\xi)|}{|t|^{p-1}}=0, \ \forall \  \xi\in\R^{N}.$
\item[$ (f_{3}) $] There exists $q\in(p,p^{*})$ such that $\lim\displaystyle_{|t|\rightarrow \infty}\dfrac{|f(t,|\xi|^{p-2}\xi)|}{|t|^{q-1}}=0, \ \forall  \ \xi\in\R^{N},$ where $p^*=\dfrac{Np}{N-p}\cdot$
\item[$ (f_{4}) $] There exists $\theta>p$ such that $$0<\theta F(t,|\xi|^{p-2}\xi)\leq tf(t,|\xi|^{p-2}\xi),$$ for all $t>0, \ \xi\in\R^{N},$ where $\displaystyle F(t,|\xi|^{p-2}\xi)=\int_{0}^{t}f(s,|\xi|^{p-2}\xi)ds.$

\item[$ (f_{5}) $] There exist positive real numbers $a$ and $b$
such that $$F(t,|\xi|^{p-2}\xi)\geq at^{\theta}-b,$$ for all $t>0, \ \xi\in \R^{N}.$

\item[$ (f_{6}) $] There exist positive constants $L_{1}$ and $L_{2}$ such that $$|f(t_{1},|\xi|^{p-2}\xi)-f(t_{2},|\xi|^{p-2}\xi)|\leq L_{1}|t_{1}-t_{2}|^{p-1}$$ 
for all $t_{1},t_{2}\in[0,\rho_{1}], \ |\xi|\leq \rho_{2},$
$$|f(t,|\xi_{1}|^{p-2}\xi_{1})-f(t,|\xi_{2}|^{p-2}\xi_{2})|\leq L_{2}|\xi_{1}-\xi_{2}|^{p-1}$$ 
for all  $t\in[0,\rho_{1}]$ and $ |\xi_{1}|,|\xi_{2}|\leq \rho_{2},$ where $\rho_{1}$, $\rho_{2}$ depend on $q,N$ and $\theta$. Moreover, $L_1$ and $L_2$ satisfy 
$\left( \dfrac{L_{2}}{C_{p}-L_{1}}\right)^{1/p-1} <1,$
 where  $C_p$ is the constant in the inequality \eqref{inn}.
 \end{itemize}
In the following, we state conditions on the potential $V:$
\begin{itemize}
\item[$ (V_{1}) $]  $V(x)\geq V_{0}>0$ for all $x\in\R^{N};$
\item[$ (V_{2}) $]  $V(x)$ is a continuous 1-periodic function, i.e., $V(x+y)=V(x), \ \forall y\in \Z^{N}$ and $\forall x\in \R^{N}.$
\end{itemize}

For further details about the periodic potential $V$, we refer to \cite{alves} and references therein.

Let  $$W=\{u\in W^{1,p}(\R^{N}):\int_{\R^{N}}(|\nabla u|^{p}+V(x)|u|^{p})dx<\infty\}.$$
 $W$ is a reflexive Banach space with the  norm
$$\lVert u\lVert=\left(\int_{\R^{N}} (|\nabla u|^{p}+V(x)|u|^{p})dx\right)^{1/p}.$$ 
Moreover, we have the continuous inclusions $W\hookrightarrow W^{1,p}(\R^{N})\hookrightarrow L^{s}(\R^{N})$  for all $s\in[p,p^{*}].$ For the details, we refer to \cite[ Lemma 2.1]{embed}.

Next, in the spirit of Figueiredo et al. \cite{figu}, we associate with \eqref{1.1}, a family of problems with no dependence on the gradient of solution.  To be precise, for every, $w \in W\cap C_{loc}^{1,\beta}(\R^{N})$ with $0<\beta<1$, we consider the  problem
\begin{equation}\label{1.2}
\begin{split}
-\Delta_{p}u+V(x)|u|^{p-2}u & =f(u,|\nabla w|^{p-2}\nabla w),  \ \ \ \\
u(x) & >0\ \   \forall x\in\R^{N}.
\end{split}
\end{equation}
Problem (\ref{1.2}) is variational in nature and 
the critical points of the functional 
$$I_{w}(u)=\frac{1}{p}\int_{\R^{N}}|\nabla u|^{p}dx+\frac{1}{p}\int_{\R^{N}}V(x)| u|^{p}dx-\int_{\R^{N}}F(u,|\nabla w|^{p-2}\nabla w)dx$$
are the weak solutions to   (\ref{1.2}).

To prove our main result, we will use the following inequality \cite{diaz}  :
\begin{equation}\label{inn}
\langle|x|^{p-2}x-|y|^{p-2}y,x-y\rangle\geq C_{p}|x-y|^{p},
\end{equation}
for all $x,y\in \R^{N},$ where $ \langle\cdot,\cdot \rangle $ is the usual inner product in $\R^{N}. $ Now, we state our main result:
\begin{theorem}\label{t1}
Suppose  that  the  conditions $ (f_{1})-(f_{6})$ and $(V_1), \, (V_2)$ are satisfied. Then, there exists a positive solution to  (\ref{1.1}).
\end{theorem}

\section{Proof of Theorem \ref{t1}}
This section deals with the proof of  Theorem \ref{t1}.
The proof  is divided in a series of lemmas. 
\begin{lemma}\label{l1}
Let $w \in W\cap C_{loc}^{1,\beta}(\R^{N})$ with $0<\beta<1$. Then there exist positive real numbers $\alpha$ and $\rho$ independent of $w$ such that
$$I_{w}(u)\geq \alpha>0, \ \ \forall u \in W\text{ such that }\|u\|=\rho.$$
\end{lemma}
\begin{proof}
From $(f_{2})$ we have, for any $\epsilon>0$, there exists $\delta_{1} >0$ such that 
\begin{equation}\label{e1}
|f(s,|\xi|^{p-2}\xi)|\leq \epsilon|s|^{p-1}, \ \ \ \forall |s|<\delta_{1}, \ \xi\in\R^{N}.
\end{equation}
From $(f_{3})$ we have, for any $\epsilon>0$, there exists $\delta_{2} >0$ such that 
\begin{equation}\label{e2}
|f(s,|\xi|^{p-2}\xi)|\leq \epsilon|s|^{q-1}, \ \ \ \forall |s|>\delta_{2}, \ \xi\in\R^{N}.
\end{equation}
By (\ref{e1}) and (\ref{e2}) we have,
\begin{equation}\label{e3}
|F(u,|\nabla w|^{p-2}\nabla w)|\leq \frac{1}{p}\epsilon|u|^{p}+\frac{1}{q}\epsilon|u|^{q}, \ \ \ \forall u \in W.
\end{equation}
Thus, 
$$I_{w}(u)=\frac{1}{p}||u||^{p}-\int_{\R^{N}}F(u,|\nabla w|^{p-2}\nabla w)dx.$$
It follows from (\ref{e3}) and embedding result that
$$I_{w}(u)\geq (\frac{1}{p}-c_{1}\epsilon)||u||^{p}-c_{2}\epsilon||u||^{q}.$$ Now choose $ \epsilon $ such that $\frac{1}{p}-c_{1}\epsilon>0$ and $\rho\leq\left( \frac{\frac{1}{p}-c_{1}\epsilon}{c_{2}\epsilon}\right) ^{\frac{1}{q-p}}.$ This completes the proof. 
\end{proof}

\begin{lemma}\label{l2}
Let $w \in W\cap C_{loc}^{1,\beta}(\R^{N})$ with $0<\beta<1$. Fix  $v_{0}\in C_{0}^{\infty}(\R^{N})$ with $||v_{0}||=1.$ Then $\exists$ $t_{0}>0,$ independent of $w,$ such that 
$$I_{w}(tv_{0})\leq 0,\ \ \forall t\geq t_{0}.$$ 
\end{lemma}
\begin{proof}
By $(f_{5})$ we get,
\begin{equation*} 
\begin{split}
I_{w}(tv_{0})& \leq \frac{t^p}{p}-\int_{\text{Supp}(v_0)}(at^{\theta}v_0^\theta-b)dx \\
 & = \frac{t^{p}}{p}-at^{\theta}\int_{\text{Supp}(v_0)}v_0^\theta dx+b|\text{Supp}(v_{0})|,
\end{split}
\end{equation*}
Since $\theta>p,$ the result follows.
\end{proof}
\begin{lemma}\label{l3}
Let conditions $(f_{1})-(f_{5})$ and $(V_1), \,(V_2)$ hold. Then, the Problem  \eqref{1.2} admits a positive solution $u_w\in W.$ 
\end{lemma}
\begin{proof}
Lemmas \ref{l1} and \ref{l2} tell us that the functional $I_{w}$ satisfies the geometric conditions of the Mountain Pass Theorem. Hence, by the version of Mountain Pass Theorem without (PS) conditions \cite{willem}, there exist a sequence $ \{u_{n}\}\subset W $ such that 
$$I_{w} (u_{n})\rightarrow c_{w} \ \ \hbox{and} \ \ I_{w}^{'}(u_{n})\rightarrow 0, \ \ \hbox{as} \ n\rightarrow \infty$$
where
$$c_{w}=\inf_{\gamma\in \varGamma}\max_{t\in [0,1]}I_{w}(\gamma(t))>0,$$
with 
$$\varGamma=\{\gamma\in C([0,1],W):\gamma(0)=0, \ \gamma(1)=t_{0}v_{0}\}$$
where $t_{0}$ and $v_{0}$ are as in Lemma \ref{l2}.\\
\indent By $(f_{4})$, we have $c\lVert u_{n}\lVert^{p}\leq c_{w}+\lVert u_{n}\lVert.$
This implies that $\{u_{n}\}$ is bounded in $W$, hence there exists its subsequence still  denoted by $\{u_{n}\}$, as
\begin{equation}
u_{n}\rightharpoonup u_{w} \ \ \hbox{in} \ \ W,
\end{equation}
\begin{equation}
u_{n}\rightarrow u_{w} \ \ \hbox{in} \ \ L_{loc}^{s} \ \ \hbox{for} \ \ p\leq s<p^{*}.
\end{equation}
On following the arguments from \cite[Proposition 4.4]{fre}, we obtain 
\begin{equation}
\frac{\partial u_{n}}{\partial x_{i}}(x)\rightarrow \frac{\partial u_{w}}{\partial x_{i}}(x) \ \ \hbox{a.e. in} \ \R^{N}.
\end{equation}
This implies,
\begin{equation}\label{2.7}
\nabla u_{n}(x)\rightarrow \nabla u_{w}(x) \ \ \hbox{a.e. in} \ \R^{N}.
\end{equation}
Using (\ref{2.7}), we get 
$$
|\nabla u_{n}|^{p-2}\nabla u_{n}\rightarrow |\nabla u_{w}|^{p-2}\nabla u_{w}\ \ \hbox{a.e. in} \ \R^{N}.$$
Since $\{|\nabla u_{n}|^{p-2}\nabla u_{n}\}$ is bounded in $L^{p/(p-1)}$, we get,
$$|\nabla u_{n}|^{p-2}\nabla u_{n}\rightharpoonup |\nabla u_{w}|^{p-2}\nabla u_{w}\ \ \hbox{ in}  \ L^{p/(p-1)}(\R^{N}).$$
By the definition of weak convergence, we have
$$\int_{\R^{N}}|\nabla u_{n}|^{p-2}\nabla u_{n}\nabla \varphi \  dx\rightarrow \int_{\R^{N}}|\nabla u_{w}|^{p-2}\nabla u_{w}\nabla \varphi \ dx \ \ \hbox{ for all } \varphi \in W.$$     
In view of Brezis-Lieb lemma \cite{h1}, we have     
$$\int_{\R^{N}}V(x)| u_{n}|^{p-2} u_{n} \varphi \  dx\rightarrow \int_{\R^{N}}V(x)| u_{w}|^{p-2} u_{w} \varphi \ dx\ \ \hbox{ for all } \varphi \in W.$$  
By the help of \cite{h1} and Lebesgue Generalized Theorem \cite{h2}, we get 
$$\int_{\R^{N}}f(u_{n},|\nabla w|^{p-2}\nabla w)\varphi \ dx\rightarrow\int_{\R^{N}}f(u_{w},|\nabla w|^{p-2}\nabla w)\varphi \ dx\ \ \hbox{ for all } \varphi \in W.$$  
Therefore, we have
\begin{equation}\label{aaaa}
\begin{split}
 I_{w}'(u_{w}) \varphi=\int_{\R^{N}}&|\nabla u_{w}|^{p-2}\nabla u_{w}\nabla \varphi  \ dx+\int_{\R^{N}}V(x)| u_{w}|^{p-2}u_{w} \varphi \ dx\\
&-\int_{\R^{N}}f(u_{w},|\nabla{w}|^{p-2}\nabla{w})\varphi \ dx=0, \hbox{ for all } \varphi \in W.
\end{split}
\end{equation}
This implies, $u$ is the weak solution of \eqref{1.2}. 

Let $u_{w}\not\equiv 0.$ Next, we  show that $u_{w}> 0$.  
By taking $\varphi=u_{w}^{-}$ in \eqref{aaaa}, we get
\begin{equation*}
\begin{split}
\int_{\R^{N}}|\nabla u_{w}|^{p-2}\nabla (u_{w}^{+}&-u_{w}^{-})\nabla u_{w}^{-}  \ dx+\int_{\R^{N}}V(x)| u_{w}|^{p-2}(u_{w}^{+}-u_{w}^{-})u_{w}^{-} \ dx\\
&=\int_{\R^{N}}f(u_{w},|\nabla{w}|^{p-2}\nabla{w})u_{w}^{-} \ dx,
\end{split}
\end{equation*}
which gives
\[
    -\int_{\R^{N}}|\nabla u_{w}|^{p-2}|\nabla u_{w}^{-}|^{2}  \ dx
        =\int_{\R^{N}}f(u_{w},|\nabla{w}|^{p-2}\nabla{w})u_{w}^{-} \ dx
        +\int_{\R^{N}}V(x)| u_{w}|^{p-2} (u_{w}^{-})^{2}  \ dx.
\]
Thus,
$$\int_{\R^{N}}V(x)| u_{w}|^{p-2} (u_{w}^{-})^{2}=0.$$
This implies, $| u_{w}|^{p-2} (u_{w}^{-})^{2}=|u_{w}^{+}-u_{w}^{-}|^{p-2} (u_{w}^{-})^{2} =0$ as $V(x)>0$. 
Therefore, we have
\[0=|u_{w}^{+}(x)-u_{w}^{-}(x)|^{p-2} (u_{w}^{-}(x))^{2} =\begin{cases} 0 \ ; &  u_{w}(x)\geq 0\\ 
|u_{w}^{-}|^{p} \ ;& u_{w}(x)< 0.
 \end{cases}\]
 Hence, $u_{w}=u_{w}^{+}-u_{w}^{-}=u_{w}^{+}\geq 0.$
Moreover,  Harnack inequality implies that $u_{w}(x)>0$ for all $x\in\R^{N}.$

If  $u_{w}\equiv 0,$ then  there exist a sequence $\{z_{n}\}\subset \R^{N}$ and $\delta, R>0$ such that
\begin{equation}\label{inter}
\int_{B_{R}(z_{n})}|u_{n}|^{p}dx\geq \delta.
\end{equation}
For, if on the contrary 
$$\limsup_{n\rightarrow\infty \ x\in \R^{N}}\int_{B_{R}(x)}|u_{n}|^{p}dx=0,$$
then by using \cite[Lemma 1.1]{Lions}, $u_{n}\rightarrow 0$ in $L^{s}(\R^{N})$ with $p<s<p^{*},$
which implies that $I_{w}(u_{n})\rightarrow 0$ as $n\rightarrow\infty.$ It contradicts the fact that $c_{w}>0$.

Let us define $v_{n}(x)=u_{n}(x+z_{n}).$ Since $V(x)$ is a 1-periodic function, we can use the invariance of $\R^{N}$ under translations to conclude that $I_{w}(v_{n})\rightarrow c_{w}$ and $I_{w}'(v_{n})\rightarrow 0$. Moreover, up to a subsequence, $v_{n}\rightharpoonup v_{w}$ in $W$ and $v_{n}\rightarrow v_{w}$ in $L^{p}(B_{R}(0)),$ where  $v_{w}$ is a critical point of $I_{w}$.  By \eqref{inter}, we conclude that $v_w$ is non zero. Arguing as above, we get that $v_w$ is a positive solution to \eqref{1.2}. This completes the proof.
\end{proof}

\begin{lemma}\label{l5}
Let  $w \in W\cap C_{loc}^{1,\beta}(\R^{N})$ with $0<\beta<1$. Then there exists positive real number $\eta$  independent of $w$, such that
$$||u_{w}||\leq \eta,$$ 
where $u_{w}$ is the solution of \eqref{1.2} obtained  in Lemma \ref{l3}.
\end{lemma}
\begin{proof}
Using the characterization of  $c_{w}$, we have
$$c_{w}\leq\max_{t\geq 0}(tu).$$
Fix $v\in W$ such that $||v||=1.$ By $(f_{5})$, we have 
$$c_{w}\leq \max_{t\geq 0}I_{w}(tv)\leq \max_{t\geq 0}\left(\frac{t^{p}}{p}-c_{6}t^{\theta}-c_{7} \right) =\eta_{0}.$$
By $(f_{4})$, we have 
\begin{equation}\label{eq4}
I_{w}(u_{w})\geq\frac{1}{p}||u_{w}||^{p}-\frac{1}{\theta}\int_{\R^{N}}f(u_{w},|\nabla w|^{p-2}\nabla w)u_{w} dx.
\end{equation}
Also, we have 
\begin{equation}\label{eq5}
 I_{w}'(u_{w})(u_{w})=||u_{w}||^{p}-\int_{\R^{N}}f(u_{w},|\nabla w|^{p-2}\nabla w)u_{w} dx.
 \end{equation}
 By (\ref{eq4}) and (\ref{eq5}), we obtain
 $$\left(\frac{1}{p}-\frac{1}{\theta}\right)||u_{w}||^{p}\leq I_{w}(u_{w})-\frac{1}{\theta}I_{w}'(u_{w})(u_{w}).$$
 Next, by using the fact that $I_{w}'(u_{w})(u_{w})=0$ and $I_{w}(u_{w})=c_{w}$ we get
  $$\left(\frac{1}{p}-\frac{1}{\theta}\right)||u_{w}||^{p}\leq c_{w}\leq \eta_{0},$$ and the proof is complete.
\end{proof}
\begin{lemma}\label{reg}
If $u_w$ is a positive solution of the equation~\eqref{1.2} obtained in Lemma~\ref{l3}, then $u_{w}\in C_{loc}^{1,\beta}\cap L_{loc}^{\infty}(\R^{N})$ with $0<\beta<1$. Moreover, there exist positive numbers $\rho_{1}$ and $\rho_{2}$, independent of $w$, such that $\lVert u_{w}\lVert _{C_{loc}^{0,\beta}(\R^{N})}\leq \rho_{1}$ and $\lVert \nabla u_{w}\lVert _{C_{loc}^{0,\beta}(\R^{N})}\leq \rho_{2}$.
\end{lemma}
\begin{proof}

By using the fact that $V(x)\geq V_{0}$ and $u_{w}> 0$, we have
$$f(u_{w},|\nabla w|^{p-2}\nabla w)-V(x)|u_{w}|^{p-2}u_{w}\leq f(u_{w},|\nabla w|^{p-2}\nabla w)-V_{0}|u_{w}|^{p-1}.$$
By the help of $(f_{2})$ and $(f_{3})$, one gets
\[
    |f(u_{w},|\nabla w|^{p-2}\nabla w)|
        \le \epsilon |u_{w}|^{p-1}+ \epsilon|u_{w}|^{q-1}+V_{0}|u_{w}|^{p-1}
        \le (\epsilon+V_{0})(|u_{w}|^{p-1}+ |u_{w}|^{q-1}).
\]
By \cite[Theorem 2.2]{regularity}, for any compact set $K\subseteq \R^N,$ we have $\| u_{w}\|_{L^{\infty}(K)}\leq C,$ where  the constant $C$  depends on $p,q,N$ and $\|u_{w}\|_{L^{p^{*}}(K)}$. By Sobolev embedding theorem and Lemma \ref{l5}, there exist $C_{0}$ independent of $w$ such that $\| u_{w}\|_{L^{\infty}(K)}\leq C_{0}$.
By \cite[Theorem~1]{ben},  $\|\nabla u_{w}\|_{L^{\infty}(K)}\leq C_{1}$, for some constant  $C_{1}$  dependent on $p,q,N$ and $\|u_{w}\|_{L^{\infty}(K)}$. Hence there exists a constant $C_2$  independent of $w$ such that  $\|\nabla u_{w}\|_{L^{\infty}(K)}\leq C_{2}.$

 By \cite[Theorem 2]{ben}, we obtain $\| u_{w}\| _{C_{loc}^{1,\beta}(\R^{N})}\leq C_{3}$, where $C_{3}$ is dependent on $p,q,N$ and $\|\nabla u_{w}\|_{L^{\infty}(K)}$. Thus there exists a positive number $\rho$ independent of $w$ such that, $\| u_{w}\| _{C_{loc}^{1,\beta}(\R^{N})}\leq \rho.$ Subsequently,  there exist positive real numbers $\rho_{1}$ and $\rho_{2}$, independent of $w$, such that $\lVert u_{w}\lVert _{C_{loc}^{0,\beta}(\R^{N})}\leq \rho_{1}$ and $\lVert \nabla u_{w}\lVert _{C_{loc}^{0,\beta}(\R^{N})}\leq \rho_{2}$. This completes the proof.
\end{proof}

\begin{lemma}\label{l4}
Let $w \in W\cap C_{loc}^{1,\beta}(\R^{N})$ with $0<\beta<1$. Then there exists positive real number $\lambda$  independent of $w$, such that
$$||u_{w}||\geq \lambda,$$ 
where $u_{w}$ is the solution of \eqref{1.2} obtained  in Lemma \ref{l3}.
\end{lemma}
\begin{proof}
Since $u_{w}$ is the weak solution of (\ref{1.2}) obtained in Lemma \ref{l3}, for all $v\in W$, we have $I_{w}'(u_{w})(v)=0$. In particular, by putting $v=u_{w}$ we get
\begin{equation*} 
\begin{split}
\int_{\R^{N}}|\nabla u_{w}|^{p}dx+\int_{\R^{N}}V(x)| u_{w}|^{p}dx&=\int_{\R^{N}}f(u_{w},|\nabla w|^{p-2}\nabla w)u_{w}dx\\
 ||u_{w}||^{p}&=\int_{\R^{N}}f(u_{w},|\nabla w|^{p-2}\nabla w)u_{w}dx.
\end{split}
\end{equation*}
By using (\ref{e1}) and (\ref{e2}) we have,
$$||u_{w}||^{p}\leq  c_{4}\epsilon||u_{w}||^{p}+c_{5}\epsilon||u_{w}||^{q}.$$
Since $q>p$, we get $||u_{w}||\geq \left(\dfrac{1-c_{4}\epsilon}{c_{5}\epsilon} \right)^{q-p}$. This completes the proof.
\end{proof}

Now, we are in position to prove Theorem \ref{t1}.

\noindent\textbf{Proof of Theorem \ref{t1}:}
Starting with an arbitrary $u_{0}\in  W\cap C_{loc}^{1,\beta}(\R^{N})$ with $0<\beta<1$, we construct a sequence $ \{u_{n}\}\subseteq  W$ as  solution of 
\begin{equation}
  -\Delta_{p}u_{n}+V(x)|u_{n}|^{p-2}u =f(u_{n},|\nabla u_{n-1}|^{p-2}\nabla u_{n-1}), \ \ \hbox{in} \ \ \R^{N}
   \tag{$P_n$}
\end{equation}
obtained  in Lemma \ref{l3}. By Lemma \ref{reg}, $\{u_n\}\subseteq W\cap C_{loc}^{1,\beta}(\R^{N})$ with $0<\beta<1$,  $\Arrowvert u_{n}\Arrowvert_{C^{0,\beta}_{loc}(\R^{N})}\leq \rho_{1}$ and $\Arrowvert \nabla u_{n}\Arrowvert_{C^{0,\beta}_{loc}(\R^{N})}\leq \rho_{2}.$ 
Since $u_{n+1}$ is the weak solution of $(P_{n+1})$, we have 
\begin{align} \label{eq6}
\int_{\R^{N}}|\nabla u_{n+1}|^{p-2}\nabla u_{n+1}\nabla \varphi \  dx
&+\int_{\R^{N}} V(x)| u_{n+1}|^{p-2}u_{n+1} \varphi \ dx \\
&=\int_{\R^{N}}f(u_{n+1},|\nabla u_{n}|^{p-2}\nabla u_{n})\varphi  \ dx,  \qquad \forall \varphi\in W.
\end{align}
Similarly, $u_{n}$ is the weak solution of $(P_{n})$, we have 
\begin{align} \label{eq7}
\int_{\R^{N}}|\nabla u_{n}|^{p-2}\nabla u_{n}\nabla \varphi  \ dx
&+\int_{\R^{N}}V(x)| u_{n}|^{p-2}u_{n} \varphi \ dx \\
&=\int_{\R^{N}}f(u_{n},|\nabla u_{n-1}|^{p-2}\nabla u_{n-1})\varphi \ dx,  \qquad \forall \varphi\in W.
\end{align}
Set $\varphi=u_{n+1}-u_{n}.$ On subtracting \eqref{eq7} from  \eqref{eq6} and by using the inequality \eqref{inn}, we get
\begin{align*}
\| u_{n+1}-u_{n} \|^{p}
& \le \frac{1}{C_{p}}\int_{\R^{N}}[f(u_{n+1},|\nabla u_{n}|^{p-2}\nabla u_{n})-f(u_{n},|\nabla u_{n}|^{p-2}\nabla u_{n})](u_{n+1}-u_{n}) dx\\
&{\quad}+\frac{1}{C_{p}}\int_{\R^{N}}[f(u_{n},|\nabla u_{n}|^{p-2}\nabla u_{n})-f(u_{n},|\nabla u_{n-1}|^{p-2}\nabla u_{n-1})](u_{n+1}-u_{n}) dx.
\end{align*}
By using $(f_{6})$, we obtain
\[
    \| u_{n+1}-u_{n} \|^{p}
        \le \frac{L_{1}}{C_{p}}\int_{\R^{N}}|u_{n+1}-u_{n}|^{p-1}(u_{n+1}-u_{n}) dx
        +\frac{L_{2}}{C_{p}}\int_{\R^{N}}|u_{n}-u_{n-1}|^{p-1}(u_{n+1}-u_{n}) dx.
\]
On simplification, we have
$$\dfrac{C_{p}-L_{1}}{C_{p}}\Arrowvert u_{n+1}-u_{n}\Arrowvert^{p}\leq\frac{L_{2}}{C_{p}}\int_{\R^{N}}|u_{n}-u_{n-1}|^{p-1}(u_{n+1}-u_{n}) dx.$$

Thanks to H$\ddot{\hbox{o}}$lder inequality, we get
$$\Arrowvert u_{n+1}-u_{n}\Arrowvert\leq\left( \dfrac{L_{2}}{C_{p}-L_{1}}\right)^{1/p-1} \Arrowvert u_{n}-u_{n-1}\Arrowvert=:d\Arrowvert u_{n}-u_{n-1}\Arrowvert,$$
where $d=\left( \dfrac{L_{2}}{C_{p}-L_{1}}\right)^{1/p-1}$. Since $d<1$,  $\{u_{n}\}$ is a Cauchy sequence in $W,$ there exists $u\in W$ such that $\{u_n\}$ converges to $u$ in $W.$ Next,  we will prove that $u$ is a solution of the Problem \eqref{1.1}. 
Since $\Arrowvert \nabla u_{n}\Arrowvert_{C^{0,\beta}_{loc}(\R^{N})}\leq \rho_{2}$,  we have
$| \ |\nabla u_{n}|^{p-2}\nabla u_{n}\nabla \varphi|\leq \rho_{2}^{p-1}|\nabla \varphi|$. 
Then, by the help of Lebesgue's Dominated Convergence Theorem, we get
$$\int_{\R^{N}}|\nabla u_{n}|^{p-2}\nabla u_{n}\nabla \varphi \  dx\rightarrow \int_{\R^{N}}|\nabla u|^{p-2}\nabla u\nabla \varphi \ dx\ \ \hbox{ for all } \varphi \in W.$$  
In view of Brezis-Lieb lemma \cite{h1}, we have     
$$\int_{\R^{N}}V(x)| u_{n}|^{p-2} u_{n} \varphi \  dx\rightarrow \int_{\R^{N}}V(x)| u|^{p-2} u \varphi \ dx\ \ \hbox{ for all } \varphi \in W.$$  
By the help of  Lebesgue Generalized Theorem \cite{h2}, we get 
$$\int_{\R^{N}}f(u_{n},|\nabla u_{n-1}|^{p-2}\nabla u_{n-1})\varphi \ dx\rightarrow\int_{\R^{N}}f(u,|\nabla u|^{p-2}\nabla u)\varphi \ dx\ \ \hbox{ for all } \varphi \in W.$$  
Therefore, as $n\rightarrow \infty,$ \eqref{eq7} implies
\begin{equation*} 
\int_{\R^{N}}|\nabla u|^{p-2}\nabla u\nabla \varphi  \ dx+\int_{\R^{N}}V(x)| u|^{p-2}u \varphi \ dx
-\int_{\R^{N}}f(u,|\nabla u|^{p-2}\nabla u)\varphi \ dx=0,
\end{equation*}
  for all  $\varphi \in W.$
This implies that $u$ is the weak solution of Problem \eqref{1.1}. By Lemma \ref{l4},  $u>0$ in $\R^{N}.$
\qed

\subsection*{Acknowledgement}
Authors would like to thank referee for his/her valuable comments and suggestions. The second author  is supported by Science and Engineering Research Board, India under the grant no. MTR/2018/000233.


\EditInfo{%
    July 17, 2020}{%
    October 27, 2020}{%
    Diana Barseghyan}

\end{paper}